\documentclass{amsart}
\usepackage[T1]{fontenc}
\usepackage{mathtools,amssymb}
\usepackage{mlmodern}
\usepackage{enumitem}

\usepackage{zref-clever}
\zcsetup{cap=true,noabbrev}

\newtheorem{theorem}{Theorem}
\newtheorem{corollary}[theorem]{Corollary}
\AddToHook{env/corollary/begin}{\zcsetup{countertype={theorem = corollary}}}
\newtheorem{lemma}[theorem]{Lemma}
\AddToHook{env/lemma/begin}{\zcsetup{countertype={theorem = lemma}}}
\newtheorem{proposition}[theorem]{Proposition}
\AddToHook{env/proposition/begin}{\zcsetup{countertype={theorem = proposition}}}
\newtheorem{question}[theorem]{Question}
\AddToHook{env/question/begin}{\zcsetup{countertype={theorem = question}}}
\newtheorem{claim}[theorem]{Claim}
\zcRefTypeSetup{claim}{Name-sg=Claim, Name-pl=Claims}
\AddToHook{env/claim/begin}{\zcsetup{countertype={theorem = claim}}}

\AddToHook{env/fact/begin}{\zcsetup{countertype={theorem = fact}}}

\theoremstyle{definition}
\newtheorem{definition}[theorem]{Definition}
\AddToHook{env/definition/begin}{\zcsetup{countertype={theorem = definition}}}

\theoremstyle{remark}
\newtheorem*{notation}{Notation}

\newcommand{\st}{\, ;\,}
\DeclarePairedDelimiter{\setlisted}{\lbrace}{\rbrace}
\DeclarePairedDelimiterX{\set}[2]{\lbrace}{\rbrace}{#1\st#2}
\DeclarePairedDelimiterX{\seq}[2]{\langle}{\rangle}{#1\st#2}
\DeclarePairedDelimiter{\card}{\lvert}{\rvert}
\DeclarePairedDelimiter{\len}{\lvert}{\rvert}

\newcommand{\funcs}[2]{\prescript{#1}{}{#2}}
\newcommand{\forces}{\Vdash}
\newcommand{\rest}{\restriction}
\renewcommand{\P}{\mathbb{P}}
\newcommand{\Q}{\mathbb{Q}}
\newcommand{\bool}{\mathbb{B}}
\newcommand{\B}{\mathbb{B}}
\newcommand{\F}{\mathcal{F}}
\newcommand{\M}{\mathbb{M}}

\newcommand{\C}{\mathbb{C}}

\DeclareMathOperator{\Add}{Add}

\newcommand{\tless}{\triangleleft}
\newcommand{\strtless}{\tless^*}

\newcommand{\addcolon}[1]{#1:}
\newlength{\casewidth}
\newlist{caselist}{enumerate}{3}
\setlist[caselist]{format=\bfseries Case~\addcolon, wide}
\setlist[caselist, 1]{label=\arabic*}
\setlist[caselist, 2]{label=\thecaselisti.\arabic*}
\setlist[caselist, 3]{label=\thecaselistii.\arabic*}

\zcRefTypeSetup{case}{Name-sg=Case, Name-pl=Cases}
\zcsetup{countertype={caselisti = case, caselistii=case, caselistiii=case}}

\begin{document}
	\title{More nonamalgamable forcing extensions}
	\author{Miha E.\ Habi\v{c}}
	\address{Proof School, 221 Main Street, San Francisco, CA 94105}
	\email{mhabic@proofschool.org}
	\urladdr{mhabic.github.io}
	
	\author{Charles Weng}
	\address{TikTok, 177 107th Ave NE, Unit 1804, Bellevue, WA 98004, USA}
	\email{charles.weng@bytedance.com}
	
	\author{Cathy Zhang}
        \address{Institute for Computational and Mathematical Engineering, Stanford University, Stanford, CA 94305, USA}
	\email{cathy.k.zhang@stanford.edu}
	
	\thanks{The authors are grateful to the anonymous referee for their careful reading of the manuscript and helpful corrections and suggestions.}
	\begin{abstract}
		We extend the results of \cite{HHKVW2019:SetTheoreticBlockchains} on nonamalgamable forcing extensions to families of posets with wide projections.
		We also use a different coding method to obtain nonamalgamable extensions by filter-based Mathias forcing.
	\end{abstract}
	
	
	\maketitle

\section{Introduction}

The \emph{generic multiverse} of a countable transitive model \(M\) of ZFC is the smallest collection of models containing \(M\) and closed under taking forcing extensions and ground models. 
Intuitively, it is the set of models one could possibly obtain from \(M\) using forcing machinery, iteratively taking extensions and/or grounds. 
Being in the same generic multiverse gives rise to an equivalence relation on the set of all countable transitive models (the only nontrivial claim is transitivity, but if \(M_2\) is reachable from \(M_1\) by a ``zig-zag'' of forcing extensions and grounds, and \(M_3\) is similarly reachable from \(M_2\), then \(M_3\) is reachable from \(M_1\) by the concatenation of the two ``zig-zags'').

One of the central themes in studying a multiverse is to understand how various forcing extensions relate and combine. A natural question is: given multiple forcing extensions \(M[G]\), can they be simultaneously \emph{amalgamated} into a single model in the multiverse containing all of them? If yes, these extensions are called \emph{amalgamable}; if no, they are \emph{nonamalgamable}.

\begin{definition}
  A family \(\mathcal{E}\) of forcing extensions of \(M\) is \emph{amalgamable} if there exists a forcing extension \(M[H]\) such that each model in \(\mathcal{E}\) is a submodel of \(M[H]\). Otherwise, \(\mathcal{E}\) is \emph{nonamalgamable}.
\end{definition}


A key tool in constructing nonamalgamable families is the use of \emph{catastrophic reals}:

\begin{definition}
  A real \(z\) is \emph{catastrophic} for \(M\) if \(z\) does not appear in any model in the generic multiverse of \(M\).
\end{definition}

By coding partial information about a catastrophic real \(z\) into each of a family of generics, we can ensure that these generics cannot all appear in a single amalgamating extension.

Catastrophic reals exist for any countable model \(M\).
One can take, for example, a real coding the ordinal height of \(M\).

There is a rich history of constructing nonamalgamable extensions. Mostowski \cite{Mostowski1976:ModelsOfBGAxioms} essentially showed that any finite poset can be embedded into the generic multiverse in a way that reflects nonamalgamability. Later, Hamkins \cite{Hamkins2016:UpwardClosureAndAmalgamationInGenericMultiverse} demonstrated that this phenomenon is pervasive across many forcing notions.

We will be particularly interested in \emph{wide posets}:

\begin{definition}
  A poset \(\P\) is \emph{wide} if, for every \(p \in \P\), there is a maximal antichain below \(p\) whose cardinality is \(|\P|\).
\end{definition}

Hamkins’ work showed that nonamalgamability patterns emerge robustly in forcing extensions by wide posets (see \cite{Hamkins2016:UpwardClosureAndAmalgamationInGenericMultiverse} for details).

A fundamental result in this direction was established in \cite{HHKVW2019:SetTheoreticBlockchains}, where the authors completely characterized which patterns of nonamalgamability can arise from families defined by \emph{finite obstacles} (see \zcref{obstaclesDef}) using wide posets.

\begin{theorem}[{\cite[Theorem 3.2]{HHKVW2019:SetTheoreticBlockchains}}]
  \label{stat:MainBlockchainThm}
  Suppose that a family of sets \(\mathcal{A}\) is defined in \(M\) by finite obstacles on a set \(I\) and \(\{\P_i : i \in I\} \in M\) is a family of wide posets in \(M\), all of the same size \(\kappa \geq |I|\). Then there are generic filters \(G_i \subseteq \P_i\) over \(M\) with the following properties:
  \begin{enumerate}
    \item If \(A \in \mathcal{A}\) then \(\prod_{i \in A} G_i\) is generic for \(\prod_{i \in A} \P_i\) over \(M\).
    \item If \(B \in M\), \(B \subseteq I\), and \(B \notin \mathcal{A}\) then the family \(\{M[G_i] : i \in B\}\) does not amalgamate in the generic multiverse.
    \item If \(A^0, A^1 \in \mathcal{A}\) then \(M[\prod_{i \in A^0} G_i] \cap M[\prod_{i \in A^1} G_i] = M[\prod_{i \in A^0 \cap A^1} G_i]\).
  \end{enumerate}
\end{theorem}


We push these results further in two significant ways:
\begin{enumerate}
  \item We show that the requirement of wide posets in \zcref{stat:MainBlockchainThm} can be relaxed. More precisely, it suffices that each poset \(\Q_i\) \emph{projects onto} a wide poset. Since many standard forcing notions (e.g., those adding Cohen reals) admit projections onto wide forcings, this shows that intricate nonamalgamability patterns occur far more broadly than previously established.
  \item We explore nonamalgamability arising from filter-based Mathias forcing. This adds a new class of examples not covered by the original wide poset framework. The techniques developed here also allow us to mix Cohen and Mathias-type extensions to witness nonamalgamability in more complex scenarios.
\end{enumerate}


In \zcref{sec:BackgroundSection} we recall the necessary background and introduce notation, focusing on forcing projections and other technical tools. 
In \zcref{sec:WideProjections} we establish that posets that project to wide posets still exhibit complex (non)amalgamability patterns. 
In \zcref{sec:MathiasSection} we present new arguments on nonamalgamability using filter-based Mathias forcing.

	\section{Background and notation}
	\label{sec:BackgroundSection}
	\subsection{Useful results and notation}
	Many constructions in this paper involve indexed collections of posets and projection maps, along with their various products.
	To avoid unnecessary notational baggage, we make the following notational convention.
	
	\begin{notation}
		Suppose \(x_i\) is a collection of sets indexed by \(i\in I\).
		Then \(x_I\) will denote the Cartesian product \(\prod_{i\in I}x_i\).
		The support of this product will be made clear in each case.%
		\footnote{We may write \(x^I\) if all \(x_i\) are the same.}
		More generally, if \(J\subseteq I\) then \(x_J\) denotes the restriction/projection of the product \(x_I\) to the coordinates in \(J\).
		
		Similarly, if \(f_i\colon x_i\to y_i\) is an indexed collection of functions and \(J\subseteq I\) then \(f_J\colon x_J\to y_J\) is the function defined on each coordinate \(j\) by \(f_j\).
	\end{notation}
	
	The following notion of strengthening of a condition in a product poset will be useful.
	
	\begin{definition}
		Let \(I\) be a set and \(\P_i\) a poset for each \(i\in I\). 
		If \(p,q\in \P_I\) and \(J\subseteq I\), we say that \(p\) is a \emph{\(J\)-extension} of \(q\), and write \(p\leq_J q\), if \(p\leq q\) and \(p\rest I\setminus J=q\rest I\setminus J\).
	\end{definition}
	
	\begin{lemma}[{\cite[Lemma 2.9]{HHKVW2019:SetTheoreticBlockchains}}]
		\label{stat:ProductDecision}
		Let $J\subseteq I$ be sets and let $\P_i$ for $i\in I$ be posets. Let $\P_I$ be the product, using whatever support you like.
		Let $\chi$ be a $\P_I$-name such that 
		$\P_I\forces \chi\subseteq V[\dot{G}_J]$ but $\P_I\not\forces\chi\in V[\dot{G}_J]$.
		Fix a condition $q\in\P_I$ such that $q\forces \chi\notin V[\dot{G}_J]$. Then there is a
		$\P_J$-name $\rho$ such that no $J$-extension of $q$ decides $\rho\in\chi$.
	\end{lemma}
	
	If \(\P\) is a separative poset, we denote by \(\bool(\P)\) its Boolean completion.
	We identify \(\P\) with a dense subset of \(\bool(\P)\) and call the conditions in \(\P\) \emph{pure} (as opposed to arbitrary conditions in \(\bool(\P)\) which are joins of antichains of pure conditions).
	By extension, a condition in a product of Boolean completions is pure if each of its coordinates is pure.
	
	We work throughout with the version of Cohen forcing \(\Add(\omega,1)=\funcs{<\omega}{2}\).
	We also fix, for the remainder of the paper, a countable transitive model \(M\) and a catastrophic real \(z\) for it.
	
	\subsection{Forcing projections}
	\label{sec:ProjectionSection}
	The idea of forcing projections appears often throughout the literature, but unfortunately there seem to be a number different (and inequivalent) definitions in use.
	Moreover, it has proven difficult to find the necessary facts spelled out (with proofs) (\cite[Section 5]{Cummings:Handbook} has a list of properties, but no proofs, while \cite{Abraham:Handbook} contains a lot of details but uses a definition stronger than ours).
	We attempt here to fix a definition for our use and collect the facts we will need. 
	
	\begin{definition}
		If \(\P\) and \(\Q\) are posets with largest elements \(1_{\P}\) and \(1_{\Q}\), respectively, we say that an order-preserving map \(\pi\colon \Q\to \P\) is a \emph{projection} if
		\begin{enumerate}
			\item \(\pi(1_{\Q})=1_{\P}\) and
			\item for any \(q\in\Q\) and any \(p\leq \pi(q)\) there is \(\bar{q}\leq q\) such that \(\pi(\bar{q})\leq p\) (in other words, the \(\pi\)-image of the cone below \(q\) is dense below \(\pi(q)\)).
		\end{enumerate}
	\end{definition}
	
	Projections come into play whenever one forcing notion adds a generic for another.

	\begin{proposition}
  If \(\pi: \Q \to \P\) is a projection and \(H \subseteq \Q\) is \(\Q\)-generic over \(M\), then \(\pi[H]\) generates a \(\P\)-generic filter over \(M\).
\end{proposition}

\begin{proof}
 Since \(H\) is a filter and \(\pi\) preserves the order, any two conditions in \(\pi[H]\) have a lower bound in \(\pi[H]\) as well.
 Moreover, it is straightforward to see that the preimage \(\pi^{-1}[D]\) of any dense open subset \(D\subseteq \P\) is dense in \(\Q\).
 It follows that the filter generated by \(\pi[H]\subseteq \P\) is generic.
\end{proof}
	
	\begin{proposition}
		\label{stat:GenericProjection}
		Let \(\Q\) and \(\P\) be posets and assume that forcing with \(\Q\) adds a generic for \(\P\). 
		More precisely, assume that there is a \(\Q\)-name \(\tau\) such that \(\Q\forces \text{\(\tau\) is \(\check{\P}\)-generic over \(\check{V}\)}\).
		Then there are a condition \(p_0\in \bool(\P)\) and a projection \(\pi\colon \Q\to \bool(\P)\rest p_0\).
	\end{proposition}
	
	When talking about \(\tau\) being generic over \(V\), the reference to the proper class ground model can be avoided by talking instead about the set of all dense subsets of \(\P\) from the ground model.
	
	\begin{proof}
  Since a generic for \(\bool(\P)\) (below a condition) can be extracted from a generic for \(\P\) (below a corresponding condition) and vice versa, we may assume at the outset that \(\P\) is a complete Boolean algebra.
  Given a condition \(q\in\Q\), let \(S_q=\set{p\in\P}{q\forces p\in\tau}\) and define a map \(\pi\colon \Q\to\P\) by letting \(\pi(q)=\bigwedge S_q\).
  Let \(p_0=\pi(1_{\Q})\).
  We claim that \(p_0\) and \(\pi\) are as desired.
  We arranged matters so that \(\pi\) maps the top condition correctly, and it is clear that it preserves the ordering as well.
  
  To check that \(\pi\) also satisfies the cone-density requirement, it is helpful to first notice that any \(q\in \Q\) will force that \(\pi(q)\in S_q\).
  This is because \(q\) forces \(S_q\) to be a (ground model) subset of the generic ultrafilter \(\tau\), which is closed under meets of ground model sets.
  So we can equivalently describe \(\pi(q)\) as the least element of \(\P\) forced by \(q\) to be in \(\tau\).
  
  With this in place, pick \(p<\pi(q)\).
  We need to find a condition \(\bar{q}\leq q\) such that \(\pi(\bar{q})\leq p\).
  If no such \(\bar{q}\) exists, then \(q\) must force \(p\notin\tau\) or, equivalently \(\lnot p\in\tau\) (where \(\lnot p\) denotes the Boolean complement of \(p\)).
  But then \(q\) also forces \(\pi(q)\wedge\lnot p\) into \(\tau\), which is a contradiction since this condition is strictly stronger than \(\pi(q)\).
\end{proof}
	
	In some cases we will be interested in, we can omit the restriction to a cone in the codomain of the projection.
	
	\begin{corollary}
		\label{stat:CohenProjection}
		Let \(\Q\) be a poset which adds a Cohen real. Then there is a projection \(\pi\colon \Q\to \bool(\Add(\omega,1))\).
	\end{corollary}
	
	\begin{proof}
  By definition, \(\Q\) adds a Cohen real, which means there is a \(\Q\)-name \(\tau\) such that \(\Q \forces \tau \text{ is } \Add(\omega,1)\text{-generic}\).
  
  Applying the \zcref{stat:GenericProjection}, there exists a condition \(p_0 \in \bool(\Add(\omega,1))\) and a projection:
  \[
    \pi': \Q \to \bool(\Add(\omega,1)) \rest p_0.
  \]
  
  The Boolean algebra \(\bool(\Add(\omega,1))\) is cone-homogeneous; that is, for any non-zero condition \(p \in \bool(\Add(\omega,1))\), there exists an isomorphism of \(\bool(\Add(\omega,1))\) and \(\bool(\Add(\omega,1))\restriction p\).
  
  Using this homogeneity, we can derive the desired projection \(\pi\) as follows:
  let \(f\colon \bool(\Add(\omega,1))\restriction p \to \bool(\Add(\omega,1))\) be an isomorphism.
  Define \(\pi = f \circ \pi'\). Then:
    \[
      \pi: \Q \to \bool(\Add(\omega,1))
    \]
    is a projection. 
    It clearly preserves the order and satisfies the two defining properties:
    \begin{enumerate}
      \item \(\pi(1_\Q) = f(\pi'(1_\Q)) = f(1_{\bool(\Add(\omega,1)) \rest p_0}) = f(p_0) = 1_{\bool(\Add(\omega,1))}\).
      \item Fix \(q \in \Q\) and \(p \leq \pi(q)\) in \(\bool(\Add(\omega,1))\).
      Write \(p' = f^{-1}(p)\).
      Note that \(p'\leq p_0\).
      Since \(\pi'\) was a projection into \(\bool(\Add(\omega,1)) \rest p_0\), there exists \(\bar{q} \leq q\) in \(\Q\) with \(\pi'(\bar{q}) \leq p'\). 
      Then:
        \[
          \pi(\bar{q}) = f(\pi'(\bar{q})) \leq f(p') = p\,.\qedhere
        \]
    \end{enumerate}  
\end{proof}
	
	\section{Posets with wide projections}	
	\label{sec:WideProjections}
	\subsection{Nonamalgamable extensions containing Cohen reals}
	\label{sec:NonamalgCohenExtensions}
	We would like to show that any posets that project to wide posets can be used in \zcref{stat:MainBlockchainThm}.
	The proof will be similar to that given in \cite{HHKVW2019:SetTheoreticBlockchains}, except that we will need to separately build generics for the original posets while also ensuring the nonamalgamability of the projected wide generics.
	
	Since the argument requires bookkeeping both in the original posets and the ones being projected to, we first sketch two easier results, where all the posets involved add a Cohen real (compare with \cite[Proposition 2.2, Theorem 2.7]{HHKVW2019:SetTheoreticBlockchains}). 
	
	The following definition will simplify some notation for us.
	
	\begin{definition}
		Let \(\pi\colon \Q\to\P\) be a projection between two posets. 
		We say that a pair \((q,p)\in \Q\times \P\) is a \emph{\(\pi\)-tagged condition} if \(p\leq \pi(q)\).
		We call \(q\) the \emph{working part} and \(p\) the \emph{tag} of the tagged condition, and omit the mention of \(\pi\) if it is clear from context.
		
		Given two tagged conditions \((q_1,p_1)\) and \((q_2,p_2)\), we say that \((q_2,p_2)\tless (q_1,p_1)\) if either:
		\begin{itemize}
			\item \(q_2=q_1\) and \(p_2\leq p_1\), or
			\item \(q_2\leq q_1\) and \(\pi(q_2)\leq p_1\).
		\end{itemize}
		We also define \((q_2,p_2)\strtless(q_1,p_1)\) if the second alternative holds.
	\end{definition}
	
	The relation \(\tless\) is a partial order on tagged conditions.
	When we refer to \emph{stronger} or \emph{weaker} tagged conditions, this is the ordering we will have in mind.
	Note that if \(q_2=q_1\) and \(\pi(q_2)\leq p_1\), it follows from the definition of a tagged condition that \(p_2\leq \pi(q_2)=\pi(q_1)=p_1\).
	That is, the second alternative implies the first one in the special case when \(q_1=q_2\).
	The idea is that we can strengthen a tagged condition by either strengthening the tag (without changing the working part), or by strengthening the working part in a way compatible with the tags.
	
	The relation \(\strtless\) fails to be reflexive (it is still antisymmetric and transitive) but is nevertheless of key interest to us.
	It is straightforward to see that, given a \(\strtless\)-descending sequence of tagged conditions, the filter generated by the working parts projects precisely to the filter generated by the tags.
	This is why tagged conditions will be useful for us: they organize a construction of a generic filter for \(\Q\) while allowing control of the projected generic for \(\P\).
	
	Unfortunately, our constructions will sometimes require us to use the weaker ordering \(\tless\).
	Still, this will not be an obstacle, provided that we use the \(\strtless\) ordering sufficiently often.
	This is implied by the following easy to check lemma.
	
	\begin{lemma}
		\label{stat:tlessTransitive}
		Suppose that \((q_1,p_1)\), \((q_2,p_2)\), and \((q_3,p_3)\) are tagged conditions.
		Then both
		\[(q_3,p_3)\tless (q_2,p_2)\strtless (q_1,p_1)\]
		and
		\[(q_3,p_3)\strtless (q_2,p_2)\tless (q_1,p_1)\]
		imply that \((q_3,p_3)\strtless (q_1,p_1)\).
	\end{lemma}
	
	
	\begin{corollary}
		\label{stat:tlessProjections}
		Suppose that \((q_n,p_n)\) form a descending sequence of tagged conditions.
		Assume also that infinitely often the strong relation \(\strtless\) appears between successive terms of this sequence.
		Then the filter generated by the working parts \(q_n\) projects to the filter generated by the tags \(p_n\).
	\end{corollary}
	
	\begin{proof}
		\zcref{stat:tlessTransitive} implies that the descending sequence of tagged conditions contains a subsequence descending in the \(\strtless\) order.
		By our earlier observation the tags of this subsequence generate the projection of the filter generated by the working parts.
		The conclusion then follows.
	\end{proof}
	
	Another important observation is that any given tagged condition \((q,p)\) can always be nontrivially strengthened, even in the strong relation \(\strtless\) (assuming that the posets involved do not have atoms).
	This is precisely because of the defining property of a projection.
	
	\begin{proposition}
		\label{stat:TwoCohenProjectionsNonamalg}
		Let \(\Q^1\) and \(\Q^2\) be forcing notions in \(M\), both of which add Cohen reals. 
		Then there are generic filters \(G^i\subseteq\Q^i\) over \(M\) such that \(M[G^1]\) and \(M[G^2]\) do not amalgamate in the generic multiverse.
	\end{proposition}
	
	We write \(\len{p}\) for the length of a Cohen condition \(p\).
	
	\begin{proof}
		Let us write \(\C=\B(\Add(\omega,1))\). 
		We call a tagged condition pure if its tag is pure.
		By \zcref{stat:CohenProjection} we have two projections \(\pi^i\colon \Q^i\to \C\). 
		We will build descending sequences of pure tagged conditions \((q^1_n,p^1_n)\in \Q^1\times \C\) and \((q^2_n,p^2_n)\in \Q^2\times \C\).
		The working parts \(q^i_n\) will generate the generics \(G^i\) and the tags will facilitate the coding of information resulting in nonamalgamability. 
		Concretely, if we think of the tags as filling in two columns side by side with 0s and 1s, we will make sure that in most rows at most one of the columns is allowed to have a 1.
		Having a row with two 1s will be a signal that the bit in the next row (in either column) is the next bit of a coded real.
		
		Let \(\seq{D^1_n}{n<\omega}\) and \(\seq{D^2_n}{n<\omega}\) be enumerations of the dense open subsets of \(\Q^1\) and \(\Q^2\) in \(M\). 
		Let \(q^i_0\) and \(p^i_0\) be the trivial conditions.
		Throughout our construction we will ensure that \(\len{p^1_n}=\len{p^2_n}\), along with the coding arrangement described above.
		So assume now that we've built the (pure) tagged conditions \((q^i_n,p^i_n)\).
		
		First, find a pure tagged condition \((q'^1,p'^1)\strtless (q^1_n,p^1_n)\) such that \(q'^1\in D^1_n\). 
		This is possible since \(\pi^1\) is a projection, \(D^1_n\) is dense open in \(\Q^1\), and the pure conditions are dense in \(\C\).
		Let \(p^*\leq p^2_n\) be the Cohen condition of length \(\len{p'^1}\) obtained by padding \(p^2_n\) with 0s.
		This makes sense since we assumed that \(\len{p^1_n}=\len{p^2_n}\) and \(\len{p'^1}\geq \len{p^1_n}\).
		Note that \((q^2_n,p^*)\) is still a pure tagged condition and \((q_n^2,p^*)\tless (q_n^2,p_n^2)\).
		We can find some pure tagged condition \((q'^2,p'^2)\strtless (q^2_n,p^*)\) and it follows that \((q'^2,p'^2)\strtless (q^2_n,p^2_n)\) as well.
		
		At this point we have \(\len{p'^2}\geq \len{p'^1}\) and we have not coded any extra information into the tags according to our scheme. 
		We now swap the roles of \(\Q^1\) and \(\Q^2\) and find pure tagged conditions \((q''^i,p''^i)\strtless (q'^i,p'^i)\) such that:
		\begin{itemize}
			\item \(q''^2\in D^2_n\) and
			\item \(p''^1\) extends the condition of length \(\len{p''^2}\) obtained by padding \(p'^1\) with 0s.
		\end{itemize}
		We can let \(q_{n+1}^i=q''^i\).
		To obtain the tags \(p_{n+1}^i\) we start with the tags \(p''^i\), pad them with 0s to the same length, then add a 1 to both of them and finally the bit \(z(n)\).
		Then \((q_{n+1}^i,p_{n+1}^i)\tless (q''^i,p''^i)\), so \((q_{n+1}^i,p_{n+1}^i)\strtless (q_n^i,p_n^i)\) by \zcref{stat:tlessProjections}.
		
		Let \(G^i\) be the filters generated by the working parts \(q_n^i\).
		Our construction ensures that the catastrophic real \(z\) can be computed from the pair of reals/filters generated by the tags \(p_n^i\) and therefore these filters cannot be found in any single model in the generic multiverse.
		However, since these filters are exactly the projections of the filters \(G^i\) by \zcref{stat:tlessProjections}, the models \(M[G^1]\) and \(M[G^2]\) do not amalgamate.
	\end{proof}
	
	The proofs in the remainder of this section will require us to construct sequences of tagged conditions in products of posets.
	In our setup, projection maps on each coordinate will give rise to a global projection map defined on the whole product poset.
	When building descending sequences of tagged conditions it will be useful for us to keep track of the coordinates on which the stronger relation \(\strtless\) holds.
	
	\begin{definition}
		Let \(\Q_i\) and \(\P_i\), for \(i\in I\), be posets and let \(\pi_i\colon \Q_i\to \P_i\) be projections.
		Let \(\pi\colon \Q_I\to \P_I\) be the induced projection.
		Let \((q_1,p_1),(q_2,p_2)\in \Q_I\times \P_I\) be \(\pi\)-tagged conditions.
		If \(J\subseteq I\) we write \((q_2,p_2)\strtless_J (q_1,p_1)\) if \((q_2,p_2)\tless (q_1,p_1)\) and \((q_2,p_2)\restriction J\strtless (q_1,p_1)\restriction J\).
	\end{definition}
		
	\begin{definition}
		\label{obstaclesDef}
		Let \(I\) be a set. A family \(\mathcal{A}\) is \emph{defined by a set of obstacles \(\mathcal{B}\)} on \(I\) if \(\mathcal{A}\) consists of those subsets of \(I\) that do not contain an element of \(\mathcal{B}\) as a subset.
	\end{definition}
	
	\begin{theorem}
		\label{stat:CohenProjectionsNonamalg}
		Suppose that a family \(\mathcal{A}\) is defined in \(M\) by a set of finite obstacles on a set \(I\) and \(\set{\Q_i}{i\in I}\) is a family of posets in \(M\), all of which add a Cohen real. 
		Then there are generic filters \(G_i\subseteq \Q_i\) over \(M\) with the following properties:
		\begin{enumerate}
			\item If \(A\in \mathcal{A}\) then \(G_A\) is generic for \(\Q_A\) over \(M\).
			\item If \(B\in M\), \(B\subseteq I\), and \(B\notin \mathcal{A}\) then the family \(\set{M[G_i]}{i\in B}\) does not amalgamate in the generic multiverse.
			\item If \(A,A'\in\mathcal{A}\) then \(M[G_A]\cap M[G_{A'}]=M[G_{A\cap A'}]\).
		\end{enumerate}
		All products above are taken with finite support.
	\end{theorem}
	
	\begin{proof}
		Our proof will again be motivated by \cite{HHKVW2019:SetTheoreticBlockchains}, but with the added bookkeeping we first saw in the proof of \zcref{stat:TwoCohenProjectionsNonamalg}.
		
		We first fix some notation.
		Let \(\mathcal{B}\) be the set of finite obstacles used to define \(\mathcal{A}\). 
		We may assume that \(\mathcal{B}\) is minimal in the sense that no two elements of \(\mathcal{B}\) are comparable via \(\subseteq\).
		Let \(\C=\bool(\Add(\omega,1))\).
		
		The projections \(\pi_i\colon \Q_i\to \C\), obtained from \zcref{stat:CohenProjection}, act coordinatewise to give rise to a projection \(\pi=\pi_I\colon \Q_I\to \C_I\).
		The plan is to recursively build a descending sequence of pure \(\pi\)-tagged conditions \((q_n,p_n)\in \Q_I\times \C_I\). 
		We will write \(q_n(i)\) or \(p_n(i)\) for the restriction of either the working parts or the tags to a single coordinate (or column) \(i\in I\).
		We will also ensure that each tag \(p_n\) is uniform, meaning that all restrictions \(p_n(i)\) have the same length (or are empty).
		
		Fundamentally, our recursion will be guided by an enumeration of all the necessary requirements. At each stage the enumeration will give us either
		\begin{enumerate}
			\item an \(A\in\mathcal{A}\) and an open dense subset \(D\subseteq \Q_A\),
			\item a \(B\in\mathcal{B}\), or
			\item a pair of sets \(A,A'\in\mathcal{A}\), together with a \(\Q_A\)-name \(\sigma\in M\) and a \(\Q_{A'}\)-name \(\tau\in M\) for a subset of \(M[\dot{G}_{A\cap A'}]\).
		\end{enumerate}
		We moreover assume that each of these objects in \(M\) is enumerated infinitely often.
		The ultimate result will be that the sequence of working parts \(q_n(i)\), for any fixed \(i\in I\), will generate a generic filter \(G_i\) for \(\Q_i\), and this generic filter will project via \(\pi\) to exactly the (generic) filter generated by the tags \(p_n(i)\).
		
		We will use a coding scheme similar to the one used in the proof of \zcref{stat:TwoCohenProjectionsNonamalg} to ensure that the real \(z\) can be decoded from the (projected) generics corresponding to the coordinates \(i\in B\) for any \(B\in\mathcal{B}\). To be precise, if we think of the tags as filling in one column of 0s and 1s for each \(i\in I\), we will make sure that the only time there is a row of 1s across all columns \(i\in B\) for some \(B\in\mathcal{B}\) is if we will use the next row to code the next bit of the catastrophic real \(z\). 
		The following definition and claim will help us show that we can carry out the work necessary for our construction without inadvertently disturbing this coding.
		
		Suppose that \(p,r\in\Add(\omega,I)\) are conditions and \(r\leq p\). 
		We shall say that \(r\) is a \emph{noncoding extension} of \(p\) if, given any \(B\in\mathcal{B}\) and any \(k\), having \(r(b,k)=1\) for all \(b\in B\) implies \(p(b,k)=1\) for all \(b\in B\) (in other words, the conditions \(q\rest B\) and \(r\rest B\) have exactly the same rows where they are constantly equal to 1).

		It is straightforward to see that the relation of being a noncoding extension is transitive.
		The following claim is a simple consequence of the definitions, but we will use it several times so we spell it out.
		
		\begin{claim}
			\label{stat:NoncodingExtensionClaim}
			Suppose \((q,p)\) is a pure tagged condition, \(A\subseteq I\), and \(D\subseteq \Q_A\) is dense open.
			Then there is a pure tagged condition \((q',p')\), with \(p'\) a uniform condition, such that:
			\begin{itemize}
				\item \((q',p')\strtless_A (q,p)\),
				\item \(q'\rest A\in D\), and
				\item \(p'\setminus p\) never takes value 1 on columns indexed by \(I\setminus A\).
			\end{itemize}
			Consequently, if \(A\in\mathcal{A}\) then \(p'\leq p\) is a noncoding extension.
		\end{claim}
		
		\begin{proof}[{Proof of {\zcref[noref]{stat:NoncodingExtensionClaim}}}]
			First, we can pad the tag \(p\) with 0s, if necessary, to make it uniform.
			Since \(\pi_A\) is a projection and \(p\rest A\leq \pi_A(q\rest A)\), we can find a condition \(\bar{q}\leq_A q\) such that \(\pi_A(\bar{q}\rest A)\leq p\rest A\).
			We then strengthen once again to \(q'\leq_A \bar{q}\) satisfying \(q'\rest A\in D\).
			To find the tag \(p'\) we pick any pure condition \(p'\rest A\leq \pi_A(q'\rest A)\) and let \(p'\rest (I\setminus A) = p\rest (I\setminus A)\).
			We also pad \(p'\) with 0s to make it uniform, if necessary.
			It is now easily seen that \((q',p')\strtless_A (q,p)\) and it is clear that we only added 1s in columns indexed by \(A\).
			
			The last bullet point implies that \(p'\leq p\) will be a noncoding extension if \(A\in\mathcal{A}\).
			This is because \(\mathcal{A}\) is closed under subsets, so it follows that a row of 1s could not have been added to any columns indexed by an element of \(\mathcal{B}\).
		\end{proof}
		
		Returning to the proof of the \zcref[noref]{stat:CohenProjectionsNonamalg}, we start the recursion by setting \(q_0\) and \(p_0\) to be the trivial conditions.
		Suppose we've constructed the pure tagged condition \((q_n,p_n)\) and consider the next requirement to be met.
		
		If we are handed an open dense subset \(D\subseteq \Q_A\) for some \(A\in\mathcal{A}\), we directly apply \zcref{stat:NoncodingExtensionClaim} to obtain \((q_{n+1},p_{n+1})\strtless_A (q_n,p_n)\) with \(p_{n+1}\leq p_n\) a noncoding extension.

		The second option for the recursive step is needing to handle an obstacle \(B\in\mathcal{B}\).
		Let \(k\) be the number of times \(B\) has occurred in the enumeration prior to now.
		We let \(q_{n+1}=q_n\) and build \(p_{n+1}\leq p_n\) by first padding all columns of \(p_n\restriction B\) with 0s to make them all have the same height\footnote{The tag \(p_n\) is already uniform, so this is only necessary if \(p_n\rest B\) has both empty and nonempty columns.}, then adding the bits 1 and \(z(k)\) to each column \(p_n(j)\) for \(j\in B\), and finally padding the whole resulting condition with 0s to make it uniform again.
		In this case \((q_{n+1},p_{n+1})\tless (q_n,p_n)\).
%
		
		In the third case we get \(A,A'\in\mathcal{A}\), a \(\Q_A\)-name \(\sigma\), and a \(\Q_{A'}\)-name \(\tau\), both for subsets of \(M[\dot{G}_{A\cap A'}]\).
		Recall that we are trying to achieve the equality \(M[G_A]\cap M[G_{A'}]=M[G_{A\cap A'}]\); we will use this step of the recursion to defeat a potential counterexample.
		We first apply \zcref{stat:NoncodingExtensionClaim} to find a pure tagged condition \((q',p')\strtless_A (q_n,p_n)\) such that \(q'\) decides whether \(\sigma\in M[G_{A\cap A'}]\) and such that \(p'\leq p_n\) is a noncoding extension.
		If \(q'\forces \sigma\in M[G_{A\cap A'}]\) we can simply let \((q_{n+1},p_{n+1})=(q',p')\) and proceed with the next step of the recursion (since the potential counterexample is no counterexample at all).
		
		If \(q'\forces \sigma\notin M[G_{A\cap A'}]\), we will attempt to find a stronger working part \(q_{n+1}\) which forces \(\sigma\neq \tau\).
		An application of \zcref{stat:ProductDecision} allows us to find a \(\Q_{A\cap A'}\)-name \(\rho\) such that no \((A\cap A')\)-extension of \(q'\rest A\) decides \(\rho\in\sigma\).
		We then apply \zcref{stat:NoncodingExtensionClaim} again to find a pure tagged condition \((q'',p'')\strtless_A (q',p')\) such that \(p''\leq p'\) is a noncoding extension and such that \(q''\) decides both whether \(\rho\in \sigma\) and \(\rho\in\tau\).
		We are going to consider two cases.
		
		\begin{caselist}
			\item Suppose that \(q''\) decides \(\rho\in\sigma\) and \(\rho\in\tau\) in opposite ways.
			In this case we can let \((q_{n+1},p_{n+1})=(q'',p'')\).
			We then get \((q_{n+1},p_{n+1})\strtless_A (q_n,p_n)\) and \(p_{n+1}\leq p_n\) is a noncoding extension.
			Moreover, \(q_{n+1}\) forces that \(\sigma\neq\tau\), since these two sets will differ on \(\rho\).
			
			\item Now suppose that \(q''\) decides \(\rho\in\sigma\) and \(\rho\in\tau\) the same way.
			Consider the condition
			\[\bar{q}=q'\rest (I\setminus A') \cup q''\rest A'\,.\]
			This condition decides \(\rho\in\tau\) the same way as \(q''\) does (since it matches \(q''\) on \(A'\)), but it does not decide \(\rho\in\sigma\).
			This is because \(\bar{q}\rest A\) is an \((A\cap A')\)-extension of \(q'\rest A\) and \(\rho\) was constructed so that no condition like that can decide \(\rho\in\sigma\).
			Let \(q_{n+1}\leq_A \bar{q}\) be some condition deciding \(\rho\in\sigma\) differently than \(\bar{q}\) decides \(\rho\in\tau\).
			Consequently \(q_{n+1}\) will force that \(\sigma\neq\tau\).
			It remains for us to find an appropriate tag \(p_{n+1}\).
			We let \(p_{n+1}\rest (I\setminus A)=p''\rest (I\setminus A)\) and find some pure \(p_{n+1}\rest A\leq \pi(q_{n+1})\rest A\), padding with 0s if necessary to make \(p_{n+1}\) uniform.
			It follows that \((q_{n+1},p_{n+1})\) is a pure tagged condition and that \((q_{n+1},p_{n+1})\strtless_A (q_n,p_n)\).
			Moreover, \(p_{n+1}\leq p_n\) is a noncoding extension, since \(p''\setminus p_n\) and \(p_{n+1}\setminus p''\) only have 1s in columns indexed by \(A\).
		\end{caselist}

		This finishes the recursive construction. Let us check that the obtained filters \(G_i\) satisfy the required properties. First, if \(A\in\mathcal{A}\), we have ensured in our construction that the filter \(G_A\) met all the dense subsets of \(\Q_A\). Therefore \(G_A\) really is generic for \(\Q_A\) over \(M\).
		
		We check the intersection property next. 
		Suppose that \(A,A'\in\mathcal{A}\). 
		We aim to see that \(M[G_A]\cap M[G_{A'}]\subseteq M[G_{A\cap A'}]\) (the other inclusion is trivially true).
		Let \(x\) be an element of this intersection.
		By \(\in\)-induction we can assume that \(x\subseteq M[G_{A\cap A'}]\).
		Fix a \(\Q_A\)-name \(\sigma\) and a \(\Q_{A'}\)-name \(\tau\) for \(x\).
		But we handled this pair of names at some stage of the recursion, at which point we guaranteed that either \(\sigma\) was forced into \(M[G_{A\cap A'}]\) or it was forced that \(\sigma\neq\tau\).
		The second alternative is clearly false in our situation, so it must be that \(x\in M[G_{A\cap A'}]\).
		
		Finally, we check that, if \(B\notin\mathcal{A}\), the extensions \(M[G_i]\) for \(i\in B\) do not amalgamate. 
		By assumption, this \(B\) must contain some obstacle \(B'\in\mathcal{B}\) and we handled this \(B'\) infinitely often in the course of our construction.
		Each time this occurred we extended the current tag \(p_n\) to add a row of 1s and the bit \(z(k)\) on the columns indexed by \(B'\), where \(k\) was the number of times \(B'\) had occurred in the construction so far.
		We claim that this allows one to compute \(z\) from the projected generics \(\pi_i[G_i]\) for \(i\in B'\).
		First, notice that \zcref{stat:tlessProjections} implies that the filters \(\pi_i[G_i]\) for any \(i\) are generated by the descending sequences of tags \(p_n(i)\).
		This is because \(\mathcal{A}\) contains all singletons, which implies that, given any \(i\in I\), there are infinitely many stages \(n\) of our construction that deal with a dense subset of \(\Q_i\).
		At that stage we will then produce \((q_{n+1},p_{n+1})\strtless_{\{i\}} (q_n,p_n)\), which suffices.
		Secondly, the reason we can do the decoding is that the components of \(z\) are encoded exactly as the bits of (any of) the generics \(\pi[G_j]\) for \(j\in B'\) which follow a coordinate \(\ell\) such that \(\pi[G_i](\ell)=1\) for all \(i\in B'\).
		This is because the only time we could have added such an index \(\ell\) was in the obstacle case of the construction, since the extensions \(p_{n+1}\leq p_n\) produced in the other two cases were noncoding.
		Moreover, since no other element of \(\mathcal{B}\) contains \(B'\), we need not worry that the bits of \(z\) coded on the coordinates of \(B'\) would be disrupted by the coding on the coordinates of some larger obstacle in \(\mathcal{B}\).
	\end{proof}
	
	\subsection{Nonamalgamable extensions with wide projections}
	Let us now move to the general version of our main theorem.
	
	\begin{theorem}
		\label{stat:mainNonamalgThm}
		Suppose that a family \(\mathcal{A}\) is defined in \(M\) by a set of finite obstacles on a set \(I\). Suppose that \(\set{\Q_i}{i\in I}\) and \(\set{\P_i}{i\in I}\) are families of posets in \(M\) such that:
		\begin{itemize}
			\item all \(\P_i\) are the same size \(\kappa\geq\card{I}\), 
			\item all \(\P_i\) are wide, and
			\item each \(\Q_i\) adds a generic for \(\P_i\).
		\end{itemize}
		Then there are generic filters \(G_i\subseteq\Q_i\) over \(M\) with the following properties:
		\begin{enumerate}
			\item\label{prop1} If \(A\in\mathcal{A}\) then \(G_A\) is generic for \(\Q_A\) over \(M\).
			\item\label{prop2} If \(B\in M\), \(B\subseteq I\), and \(B\notin \mathcal{A}\) then the family \(\set{M[G_i]}{i\in B}\) does not amalgamate in the generic multiverse.
			\item\label{prop3} If \(A^0,A^1\in\mathcal{A}\) then \(M[G_{A^0}]\cap M[G_{A^1}] = M[G_{A^0\cap A^1}]\).
		\end{enumerate}
		The products above may use supports from an arbitrary (fixed) ideal in \(M\) extending the finite ideal on \(I\).
	\end{theorem}
	
	
	We should mention that the heart of the above theorem (properties \zcref[noname]{prop1,prop2}) can be obtained very easily from known results.
	As was essentially observed by Hamkins \cite[Theorem 8]{Hamkins2016:UpwardClosureAndAmalgamationInGenericMultiverse}, if we can achieve the desired pattern of (non)amalgamability for the \(\P_i\)-extensions (using results from \cite{HHKVW2019:SetTheoreticBlockchains}), we simply force over those by the relevant quotient forcings to obtain \(\Q_i\)-extensions exhibiting the same pattern.
	However, more work is necessary to make sure that the quotient forcing does not break property (\zcref[noname]{prop3}).
	
	\begin{proof}
		By \zcref{stat:GenericProjection} we have, for each \(i\in I\), a projection \(\pi_i\colon \Q_i\to \B_i\), where \(\B_i\) is the Boolean completion of \(\P_i\).
		The coordinatewise projections give rise to a projection \(\pi=\pi_I\colon \Q_I\to \B_I\). 
		
		Since all \(\P_i\) have the same size \(\kappa\), we can fix enumerations of each of them in that order type within \(M\).
		Moreover, since all \(\P_i\) are wide, we can also find, for each \(p\in\P_i\) a maximal antichain \(W_i(p)\) below \(p\), together with its enumeration in order type \(\kappa\).
		
		We will build a descending sequence of pure \(\pi\)-tagged conditions \((q_n,p_n)\in \Q_I\times \B_I\), whose working parts will grow into the desired generics \(G_i\) and whose tags will ensure the nonamalgamability of the forcing extensions.
		We write \(q_n(i)\) and \(p_n(i)\) for the \(i\)th component of each of these conditions.
		At each step of the construction we will deal with: two sets \(A^0,A^1\in\mathcal{A}\), a dense subset \(D\subseteq \Q_{A^0}\) in \(M\), and a \(\Q_{A^0}\)-name \(\sigma\) and a \(\Q_{A^1}\)-name \(\tau\) in \(M\) for subsets of \(M[\dot{G}_{A^0\cap A^1}]\).%
		\footnote{It turns out to be more convenient to deal with all requirements at once rather than to split them up into different types as in the proof of \zcref{stat:CohenProjectionsNonamalg}.}
		Into this process we will also interleave the construction of an auxiliary sequence of conditions \(c_n\in \P_I\), to be used as \emph{coding points}.
		
		The proof strategy will be essentially that of \cite[Theorem 3.2]{HHKVW2019:SetTheoreticBlockchains}, with ideas from \zcref{sec:NonamalgCohenExtensions}.
		We will code information into where a particular (projected) generic containing a condition \(p\) meets the maximal antichain \(W_i(p)\).
		This amounts to coding an ordinal below \(\kappa\).
		However, using some predetermined scheme, we will in fact carry with us finitely many finite subsets of \(\kappa\).
		We will also reserve three special symbols \(\sharp\), \(\flat\), and \(\natural\), to serve as markers.
		Along the way it will be important for us to have some way of knowing which condition \(p\) to consult regarding the next step of coding.
		This is where the sequence of coding points will come in useful.
		
		Start the construction by letting \(q_0\) and \(p_0\) be the trivial conditions.
		Now suppose that we have constructed the pure tagged condition \((q_n,p_n)\) and are faced with \(A^0_n, A^1_n\in\mathcal{A}\), a dense \(D\subseteq \Q_{A^0_n}\), and names \(\sigma,\tau\).
		Let \(B^0_n\) be the set of the first \(n\) elements of \(I\setminus A^0_n\) in some fixed enumeration of \(I\) in order type \(\omega\), and similarly for \(B^1_n\) (since \(I\in M\) and \(M\) is countable, there is such an enumeration in \(V\)).
		
		Lest we forget, we are also constructing the sequence of coding points \(c_n\in \P_I\).
		We assume that the construction so far satisfies the following property:
		for each \(m<n\) we have \(p_{m+1}\leq_{B^0_m} c_m\) and for each \(i\in B^0_m\) the condition \(p_{m+1}(i)\leq c_{m}(i)\) is that element of the antichain \(W_i(c_{m}(i))\) which codes
		\begin{itemize}
			\item the sets \(B^0_{m+1}\) and \(B^1_{m+1}\),
			\item the bit \(z(m)\), and
			\item the condition \(c_m\rest B^0_{m+1}\).
		\end{itemize}
		Thus, once we've built \(c_m\), the condition \(p_{m+1}\) will be completely determined (we should also take care to construct \(c_m\) so that \((q_{m+1},p_{m+1})\) is a tagged condition).
		
		We now aim to build the next coding point \(c_n\) and tagged condition \((q_{n+1},p_{n+1})\).
		First, find some \(\bar{q}\leq_{A_n^0}q_n\) such that \(\pi(\bar{q})\rest A_n^0\leq p_n\restriction A_n^0\) and \(\bar{q}\) is a condition in \(D\) deciding whether \(\sigma\in M[\dot{G}_{A^0_n\cap A^1_n}]\).
		Such a \(\bar{q}\) exists because we assumed that \(p_n\leq \pi(q_n)\) and \(\pi\) is a projection.
		Now consider a handful of cases:				
		
		\begin{caselist}
			\item \label{caseEasy}
			Suppose that \(\bar{q}\forces \sigma\in M[\dot{G}_{A^0_n\cap A^1_n}]\).
			In this case we let \(q_{n+1}=\bar{q}\) and find \(c_n\) as follows:
			let \(\bar{p}\leq_{A^0_n} p_n\) be a pure condition such that \(\bar{p}\rest A^0_n \leq \pi(q_{n+1})\rest A^0_n\).
			It follows that \((q_{n+1},\bar{p})\) is a tagged condition and \((q_{n+1},\bar{p})\strtless_{A_n^0} (q_n,p_n)\).
			We construct the coding point \(c_n\leq_{B_n^0} \bar{p}\) by letting \(c_n(i)\leq \bar{p}(i)\) be the element of \(W_i(\bar{p}(i))\) coding the symbol \(\sharp\), for \(i\in B_n^0\).
			Now build \(p_{n+1}\) from \(c_n\) as stipulated above.
			This gives us a pure tagged condition \((q_{n+1},p_{n+1})\tless (q_{n+1},\bar{p})\strtless_{A_n^0} (q_n,p_n)\).
			
			
			\item Suppose that \(\bar{q}\forces \sigma\notin M[\dot{G}_{A_n^0\cap A_n^1}]\).
			In this case we will attempt to guide our construction to force \(\sigma\neq \tau\).
			By \zcref{stat:ProductDecision} there is a \(\Q_{A_n^0\cap A_n^1}\)-name \(\rho\) such that no \((A_n^0\cap A_n^1)\)-extension of \(\bar{q}\) decides whether \(\rho\in \sigma\).
			We consider two subcases:
			
			\begin{caselist}
				\item \label{caseStillEasy}
				Suppose that \(\bar{q}\) decides whether \(\rho\in\tau\).
				Let \(q_{n+1}\leq_{A_n^0} \bar{q}\) be any extension which decides \(\rho\in\sigma\) in the opposite way and build \(c_n\) and \(p_{n+1}\) as in \zcref{caseEasy}.
				This again gives us \((q_{n+1},p_{n+1})\strtless_{A_n^0} (q_n,p_n)\).
				
				\item \label{caseHarder}
				Suppose that \(\bar{q}\) does not decide \(\rho\in\tau\).
				We first find \(q'\leq_{A^0_n} \bar{q}\) forcing that \(\rho\in\sigma\).
				As in \zcref{caseEasy} let \(p'\leq_{A_n^0} p_n\) be a pure condition such that \(p'\rest A_n^0\leq \pi(q')\rest A_n^0\).
				Then \((q',p')\) is a tagged condition and \((q',p')\strtless_{A_n^0} (q_n,p_n)\).
				Now build \(p''\leq_{B^0_n} p'\) by letting \(p''(i)\), for \(i\in B^0_n\), be the element of \(W_i(p'(i))\) coding the symbol \(\flat\) and the sequence of conditions \(p'\rest B_n^1\).
				We would like to use \(p''\) as the next coding point, but we must consider whether committing to this would still allow us to force \(\rho\notin \tau\).
				
				\begin{caselist}
					\item \label{caseHardButOk}
					Suppose that there is some condition \(q^*\leq_{A_n^1} q'\) forcing \(\rho\notin\tau\) and satisfying \(\pi(q^*)\rest A_n^1\leq p''\rest A_n^1\).
					Let \(p^*\leq \pi(q^*)\) be a pure condition such that \(p^*\leq_{A_n^1} p''\).
					Such \(p^*\) exist since \(p''\leq \pi(q')\) and we assumed that strengthening \(q'\) to \(q^*\) only involved coordinates in \(A_n^1\).
					Then \((q^*,p^*)\) is a tagged condition and \((q^*,p^*)\strtless_{A_n^1} (q',p'')\strtless_{A_n^0}(q_n,p_n)\).
					Let \(q_{n+1}=q^*\) and build \(c_n\leq_{B_n^1} p^*\) by letting \(c_n(i)\in W_i(p^*(i))\), for \(i\in B_n^1\) be the element coding the sequence of conditions \(p^*\rest B_n^0\).
					If we then build \(p_{n+1}\leq c_n\) it follows that \((q_{n+1},p_{n+1})\tless (q^*,p^*)\).
					
					\item \label{caseHardest}
					Suppose that no condition \(q^*\) exists as in \zcref{caseHardButOk}.
					In other words, any condition \(q^*\rest A_n^1\leq q'\rest A_n^1\) satisfying \(\pi(q^*)\rest A_n^1\leq p''\rest A_n^1\) forces \(\rho\in\tau\). 
					Define
					\[q^{**}=\bar{q}\restriction(A^0_n\setminus A^1_n) \cup q'\rest (I\setminus(A^0_n\setminus A^1_n))\]
					and
					\[p^{**}=\pi(\bar{q})\restriction(A^0_n\setminus A^1_n)\cup p''\rest (I\setminus(A^0_n\setminus A^1_n))\,.\]
					Notice that \((q^{**},p^{**})\) is a tagged condition and \((q^{**},p^{**})\strtless_{A^0_n} (q_n,p_n)\).
					Moreover, observe that \(q^{**}\) does not decide \(\rho\in\sigma\).
					This is because \(q^{**}\rest A^0_n\) is an \((A_n^0\cap A_n^1)\)-extension of \(\bar{q}\restriction A_n^0\) and, by our construction of \(\rho\), no such extension can decide \(\rho\in\sigma\).
					Now we find \(q_{n+1}\leq_{A_n^0} q^{**}\) forcing \(\rho\notin \sigma\) and satisfying \(\pi(q_{n+1})\rest A_n^0\leq p^{**}\rest A_n^0\).
					We do this in two steps: first we strengthen \(q^{**}\rest (A_n^0\cap A_n^1)\) to a condition that projects below \(p^{**}\rest (A_n^0\cap A_n^1)\).
					The resulting condition still does not decide \(\rho\in\sigma\) for the same reason as above.
					Secondly, we \(A_n^0\)-extend this intermediate condition to obtain \(q_{n+1}\) forcing \(\rho\notin\sigma\).
					Let \(p'''\leq_{A_n^0} p''\) be a pure condition such that \(p'''\leq \pi(q_{n+1})\).
					This gives us a tagged condition \((q_{n+1},p''')\strtless_{A^0_n} (q^{**},p^{**})\).
					Now build \(c_{n}\leq_{B_n^0} p'''\) by letting \(c_n(i)\in W_i(p'''(i))\), for \(i\in B_n^0\), be the element coding the symbol \(\natural\).
					If we then build \(p_{n+1}\leq c_n\) as before, we get \((q_{n+1},p_{n+1})\strtless_{A_n^0} (q_n,p_n)\).
%
%
%
%
%
				\end{caselist}
			\end{caselist}
		\end{caselist}
		
		This finishes the construction of the coding point \(c_n\) and the tagged condition \((q_{n+1},p_{n+1})\).
		Let \(G_i\) for \(i\in I\) be the filter generated by the descending sequence \(q_n(i)\).
		We claim that these satisfy the properties listed in the theorem.
		
		Properties (\zcref[noname]{prop1}) and (\zcref[noname]{prop3}) are easy to see.
		If \(A\in\mathcal{A}\) then we made sure that the descending sequence of restricted working parts \(q_n\rest A\) met every dense subset \(D\subseteq \Q_{A}\) in \(M\), meaning that \(G_A\) is generic for \(\Q_A\), as desired.
		
		For property (\zcref[noname]{prop3}), suppose that \(A^0,A^1\in\mathcal{A}\) and let \(x\in M[G_{A^0}]\cap M[G_{A^1}]\). 
		We wish to see that \(x\in M[G_{A^0\cap A^1}]\) and can assume, by induction, that \(x\subseteq M[G_{A^0\cap A^1}]\).
		We can also find \(\Q_{A^0}\) and \(\Q_{A^1}\)-names \(\sigma\) and \(\tau\) for \(x\).
		But these names were considered at some stage of our construction and we ensured at that stage that either \(\sigma\) was forced into \(M[G_{A^0\cap A^1}]\) (in \zcref{caseEasy}), or it was forced that \(\sigma\neq\tau\) (in the other cases).%
		\footnote{While we didn't explicitly force \(\rho\in\tau\) in \zcref{caseHardest}, we did make sure that any condition \(\strtless_{A^1_n}\)-below \((q^{**},p^{**})\) (or \((q_{n+1},p_{n+1})\)) forces that statement. This is enough, since at each step \(k\) we get \((q_{k+1},p_{k+1})\strtless_{A^0_k} (q_k,p_k)\) and we can ensure in our enumeration that any given \(A\) appears as \(A_k^0\) for infinitely many \(k\).}
		But if \(\sigma\) and \(\tau\) were both names for \(x\), it must be the case that \(x\in M[G_{A^0\cap A^1}]\) as desired.
		
		We turn now to verifying property (\zcref[noname]{prop2}). 
		So suppose that \(B\notin\mathcal{A}\) is a subset of \(I\).
		We wish to see that \(\set{M[G_i]}{i\in B}\) does not amalgamate.
		Since we assumed that \(\mathcal{A}\) is defined by finite obstacles, we may assume that \(B\) is finite (otherwise we repeat this argument for a finite obstacle contained in \(B\)).
		We will show that we can reconstruct the catastrophic real \(z\) from the sequence of projected generics \(\pi[G_B]\) (and some other parameters in \(M\)), essentially by retracing the construction of the tags \(p_n\) and the coding points \(c_n\).
		For this to work we need to make sure that the projected generics have something to do with the tags.
		This will be the case since \(\mathcal{A}\) contains all singletons and since \((q_{n+1},p_{n+1})\strtless_{A^0_n} (q_n,p_n)\) for all \(n\). \zcref{stat:tlessProjections} then implies that each projection \(\pi_i[G_i]\) is generated by the sequence of tags \(p_n(i)\).
		
		So let's assume that we know \(B_{n}^0\), \(B_{n}^1\), and \(z\rest n\), as well as the partial tag \(p_n\rest B_{n}^0\).
		We may also assume that \(n\) is large enough so that \(B\subseteq A_m^0\cup B_m^0\) and \(B\subseteq A_m^1\cup B_m^1\) for all \(m\geq n\).
		This is possible since \(B\) is finite.
		Moreover, we know that \(B\cap B_m^0\) and \(B\cap B_m^1\) are nonempty for all such \(m\), since \(\mathcal{A}\) is closed under subsets.
	
		
		According to our construction, we can pick some \(i\in B\cap B_n^0\) and check where the projected generic \(\pi_i[G_i]\) meets the antichain \(W_i(p_n(i))\). 
		If we see \(\sharp\), the construction passed through either \zcref{caseEasy} or \zcref{caseStillEasy}.
		On the other hand, if we see \(\flat\), the construction passed through \zcref{caseHarder}.
		In the first two cases, we can first recover \(c_n\rest B_n^0\) and then use the generic \(\pi_i[G_i]\) to see what is coded in the antichain below \(c_n(i)\).
		This reveals \(B_{n+1}^0\), \(B_{n+1}^1\), and \(z(n)\), as well as \(c_n\rest B^0_{n+1}\).
		We then use this to build \(p_{n+1}\rest B^0_{n+1}\), coding the same information below \(c_n\rest (B^0_{n+1}\cap B^0_n)\) and keeping the conditions \(c_n\rest (B^0_{n+1}\setminus B^0_n)\).
		
		Now suppose we've discovered the construction passed through \zcref{caseHarder}.
		In this case, we can also reconstruct \(p'\rest B_n^0\cup B_n^1\) and \(p''\rest B_n^0\cup B_n^1\) (where \(p'\) and \(p''\) are the conditions as described above).
		We now check where the projected generic \(\pi_i[G_i]\) meets the antichain \(W_i(p''(i))\).
		If we see \(\natural\), the construction passed through \zcref{caseHardest}, and through \zcref{caseHardButOk} otherwise.
		In the first case we first recover \(c_n\rest B_n^0\) and then proceed to build \(p_{n+1}\rest B^0_{n+1}\) as described above for \zcref{caseEasy,caseStillEasy}.
		In the second case we pick some \(j\in B\cap B_n^1\) and use \(\pi_j[G_j]\) to probe the antichain \(W_j(p''(j))\).
		This will reveal the sequence \(p^*\rest B_n^0\) which we can use to reconstruct \(c_n\rest B_n^0\).
		After this we use \(\pi_i[G_i]\) to probe \(W_i(c_n(i))\) and build \(p_{n+1}\rest B^0_{n+1}\) as before.
	\end{proof}
	
	\zcref{stat:mainNonamalgThm} improves on the results of \cite{HHKVW2019:SetTheoreticBlockchains} particularly in the case of uncountable ccc forcings which add Cohen reals.
	We were originally motivated to investigate the topic of this paper after proving a series of ad hoc versions of the main theorem in specific cases: Hechler forcing, the standard forcing to specialize an Aronszajn tree, etc.
	All of these forcings added Cohen reals, and are therefore covered by \zcref{stat:CohenProjectionsNonamalg}.
	We add to this collection in the next section, where we will examine the nonamalgamability of forcing extensions by a version of Mathias forcing.
	
	In general, the nonamalgamability phenomenon for ccc forcing depends on the specific forcing used.
	Jensen's forcing to add a unique generic real over \(L\) has the property that any finitely many generics for it are mutually generic and therefore does not give rise to nonamalgamable extensions.\footnote{For two generics this is \cite[Corollary 28.6]{Jech2002:SetTheory}, but the argument generalizes to an arbitrary finite collection.
	See \cite{FriedmanGitmanKanovei2019:ModelOfSecondOrderArithmetic} for an exposition of products of Jensen forcing.}
	Fuchs and Hamkins \cite{FuchsHamkins2009:DegreesOfRigiditySuslinTrees} constructed a Suslin tree with the same property using \(\diamondsuit\) (and also showed that a generic Suslin tree works as well).
	Nevertheless, it remains open whether the nonamalgamability phenomenon appears with some basic ccc forcing notions, such as random real forcing.
	
	After circulating a version of this paper, we received the following observation from Jason Zesheng Chen, constructing a pair of nonamalgamable random reals.
	
	\begin{proposition}[Chen]
		There is a pair of nonamalgamable random reals over \(M\).
	\end{proposition}
	
	\begin{proof}
		By a result of Solovay (see \cite[Theorem 3.1.3]{BartoszynskiJudah:book}) the set \(R(M)\) of random reals over \(M\) has measure 1.
		On the other hand, the set \(\set{x\oplus z}{x\in R(M)}\) is just a shift of \(R(M)\) and therefore also has measure 1.
		This means that these two sets must intersect, so there is a random real \(x\) such that \(x\oplus z\) is also random over \(M\).
		But we can recover \(z\) from \(x\) and \(x\oplus z\), so \(M[x]\) and \(M[x\oplus z]\) do not amalgamate.
	\end{proof}
	
	It is not clear whether this kind of soft argument extends to provide more complicated arrangements of nonamalgamable random extensions, so we ask the following question.
	
	\begin{question}
		Are there three random real forcing extensions that are pairwise amalgamable, but not fully amalgamable?
	\end{question}
	
	\section{Nonamalgamable Mathias-like extensions}
	\label{sec:MathiasSection}
	In this section we would like to present another example of nonamalgamability in the generic multiverse that does not (in general) follow from already established results.
	We will be concerned with a ccc version of Mathias forcing.
	The ordinary Mathias forcing is wide and therefore exhibits the full array of nonamalgamability by \zcref{stat:MainBlockchainThm} or \zcref{stat:mainNonamalgThm}.
	However, a common variant of the forcing restricts the upper parts of the conditions to a particular filter (often an ultrafilter) on \(\omega\).
	
	\begin{definition}
		Let \(\F\) be a filter on \(\omega\).
		The \emph{\(\F\)-based Mathias forcing} \(\M_\F\) has conditions \((s,A)\) where \(s\) is a finite subset of \(\omega\) and \(A\in \F\) with \(\min A>\max s\).
		The conditions are ordered by letting \((t,B)\leq (s,A)\) if \(t\) end-extends \(s\), \(B\subseteq A\), and \(t\setminus s\subseteq A\).
	\end{definition}
	
	The poset \(\M_\F\) is not interesting unless \(\F\) extends the cofinite filter, so we assume this to be the case throughout.
	The forcing \(\M_\F\) is clearly ccc (and even \(\sigma\)-centred): any two conditions \((s,A)\) and \((s,B)\) with the same stems are compatible since we can just intersect \(A\) and \(B\) and get another element of \(\F\).
	It is known that, in the case that \(\F\) is an ultrafilter, \(\M_\F\) adds a Cohen real if and only if \(\F\) is not Ramsey (see \cite{JudahShelah1991:ForcingMinimalDegree}).
	
	\begin{proposition}
		\label{stat:TwoMathiasNonamalg}
		There are generic filters \(G,H\subseteq \M_\F\) over \(M\) such that \(M[G]\) and \(M[H]\) do not amalgamate in the generic multiverse.
	\end{proposition}
	
	\begin{proof}
		We will build two descending sequences \(p_n=(s_n,A_n)\) and \(q_n=(t_n,B_n)\) of Mathias conditions.
		At each step of the construction we will meet a dense open set and code some information about \(z\) into the conditions we build.
		The coding is based on the notion of \emph{oscillation} due to Todorčević (see, for example, \cite{Todorcevic1988:OscillationsOfReals}).
		The coding signal for us will be a point \(x\) in the intersection of the two Mathias reals we are building.
		We will code a bit of the catastrophic real \(z\) by choosing which of the two Mathias reals will get the first element following a coding point.
		
		More precisely, let \(\seq{D_n}{n<\omega}\) enumerate the dense open subsets of \(\M_\F\) in \(M\).
		Let \(p_0\) and \(q_0\) be the trivial conditions.
		Now assume we've built \(p_n\) and \(q_n\).
		We first introduce a coding point by extending the stems of \(p_n\) and \(q_n\) by \(m=\min (A_n\cap B_n)\).
		Let \(p_n'=(s_n\cup\{m\},A_n\setminus (m+1))\) and similarly for \(q_n'\).
		Now consider the bit \(z(n)\).
		Suppose that \(z(n)=0\).
		In this case let \(p_{n+1}\) be any extension of \(p_n'\) in \(D_n\) which extends the stem of \(p_n'\) by at least one point.
		Then, if \(q_n'=(t_n',B_n')\), let \(q_{n+1}\) be any extension of \((t_n',B_n'\setminus (\max s_{n+1}+1))\) in \(D_n\) which extends the stem by at least one point.		
		On the other hand, if \(z(n)=1\) we swap the roles of \(p_n\) and \(q_n\), extending \(q_n\) first and \(p_n\) second.
		
		This finishes the construction of the conditions \(p_{n+1}\) and \(q_{n+1}\).
		It is clear that these descending sequences of conditions will generate generic filters \(G\) and \(H\).
		Meanwhile, from \(G\) and \(H\) (or rather, their corresponding Mathias reals) we can reconstruct the catastrophic real \(z\).
		Specifically, to recover the bit \(z(n)\) we find the \(n\)th element of the intersection of the two generic reals and then inspect the following element of the two reals.
		If the next element of the real corresponding to \(G\) is smaller than the one corresponding to \(H\) then \(z(n)=0\) and otherwise \(z(n)=1\).
		This works since our construction ensured that the only points in the intersection of the two generic reals are the coding points.
	\end{proof}
	
	We can extend this to a version of \cite[Theorem 2.7]{HHKVW2019:SetTheoreticBlockchains}.
	
	\begin{theorem}
		\label{stat:MathiasNonamalg}
		Suppose that a family \(\mathcal{A}\) is defined in \(M\) by a set of finite obstacles on a set \(I\).
		Let \(\F\in M\) be a filter on \(\omega\) extending the cofinite filter.
		Then there are generic filters \(G_i\subseteq \M_\F\) over \(M\), for \(i\in I\), with the following properties:
		\begin{enumerate}
			\item If \(A\in\mathcal{A}\) then \(G_A\) is generic for \(\M_\F^A\) over \(M\).
			\item If \(B\in M\), \(B\subseteq I\), and \(B\notin\mathcal{A}\) then the family \(\set{M[G_i]}{i\in B}\) does not amalgamate in the generic multiverse.
			\item If \(A,A'\in\mathcal{A}\) then \(M[G_A]\cap M[G_{A'}]=M[G_{A\cap A'}]\).
		\end{enumerate}
		All products above are taken with finite support.
	\end{theorem}
	
	\begin{proof}
		Let \(\mathcal{B}\) be the set of finite obstacles defining \(\mathcal{A}\).
		We may assume that \(\mathcal{B}\) consists of minimal obstacles, in the sense that no two elements of \(\mathcal{B}\) are comparable via \(\subseteq\).
		We will build a descending sequence of conditions \(p_n=(s_n,X_n)\in \M_\F^I\) by recursion, meeting dense sets and coding information about \(z\). 
		
		At each step of our construction we will be given either
		\begin{enumerate}
			\item an \(A\in\mathcal{A}\) and an open dense subset \(D\subseteq \M_\F^A\),
			\item a \(B\in\mathcal{B}\), or
			\item a pair of sets \(A,A'\in\mathcal{A}\), together with a \(\M_\F^A\)-name \(\sigma\in M\) and a \(\M_\F^{A'}\)-name \(\tau\in M\) for a subset of \(M[\dot{G}_{A\cap A'}]\).
		\end{enumerate}
		
		We will make sure to precisely control the intersection of the various \(s_n(i)\), as those will serve as our coding points.
		The coding scheme will be like the one used in \zcref{stat:TwoMathiasNonamalg}.
		We fix in advance, for each \(B\in\mathcal{B}\), an element \(i_B\in B\).
		This is the index we will consult after a coding point is reached to find the next bit of the real \(z\).
		After each step of the construction we will also make sure that the built condition is \emph{uniform}, meaning that all nontrivial coordinates have the same upper part.
		We can do this by simply intersecting the upper parts of those finitely many coordinates, since \(\F\) is a filter.
		
		We start the construction by letting \(p_0\) be the trivial condition.
		Now assume we've built \(p_n\) and consider the next task.
		If we receive a dense subset \(D\) of some \(\M_\F^A\), we simply extend \(p_n\rest A\) to meet \(D\) and make the condition uniform.
		
		If we are given a finite obstacle \(B\), let \(k\) be the number of times \(B\) has appeared in the enumeration before this.
		We first extend \(p_n\rest B\) by adding \(\min\bigcap_{i\in B} X_n(i)\) to each of the stems and then consider the value of \(z(k)\).
		If \(z(k)=0\) we extend the stem of \(p_n(i_B)\) by one point \(b\) and the stems of \(p_n\rest (B\setminus \setlisted{i_B})\) by one point each, all strictly above \(b\).
		If, on the other hand, \(z(k)=1\), we extend the stems of \(p_n\rest (B\setminus \setlisted{i_B})\) by one point each and then extend the stem of \(p_n(i_B)\) by a point above all of the (finitely many) ones just added.
		In either case, we finish by making the condition uniform.
		
		Finally, if given two names \(\sigma\) and \(\tau\) as above, we will attempt to prevent them from naming the same set not in \(M[G_{A\cap A'}]\).
		First consider whether \(p_n\rest A\forces \sigma \in M[\dot{G}_{A\cap A'}]\).
		If so, we do nothing and move on to the next step of the construction.
		Otherwise we will build \(p_{n+1}\leq p_n\) forcing \(\sigma\neq \tau\).
		By our assumption we can find a \(q\leq_A p\) such that \(q\rest A\forces \sigma\notin M[\dot{G}_{A\cap A'}]\).
		Using \zcref{stat:ProductDecision} we find a name \(\rho\) such that no \((A\cap A')\)-extension of \(q\) decides \(\rho\in\sigma\).
		Find an \(A\)-extension of \(q\) deciding \(\rho\in\sigma\) and make it uniform; call the resulting condition \(q'\).
		Now consider two cases:
		\begin{enumerate}
			\item If \(q'\rest A'\) does not decide \(\rho\in\tau\) the same way as \(q'\rest A\) decides \(\rho\in\sigma\), we can find a \(q''\leq_{A'} q'\) deciding it the opposite way.
			We then obtain \(p_{n+1}\) by making this \(q''\) uniform.
			
			\item Suppose that \(q'\rest A'\) decides \(\rho\in\tau\) the same way as \(q'\rest A\) decides \(\rho\in \sigma\). 
			It follows from our construction of \(\rho\) that \(q\rest (A\setminus A') \cup q'\rest (I\setminus (A\setminus A'))\) does not decide \(\rho\in\sigma\) but still decides \(\rho\in\tau\).
			We can thus find an \(A\)-extension \(q''\) of this condition which decides \(\rho\in\sigma\) in the opposite way.
			We then obtain \(p_{n+1}\) by making this \(q''\) uniform.
		\end{enumerate}
		
		This finishes the construction of the conditions \(p_n\).
		Let \(G_i\) be the filters generated by this sequence and \(r_i\) the corresponding reals, obtained as the unions of the stems.
		Our construction makes it clear that \(G_A\) is generic whenever \(A\in \mathcal{A}\).
		Moreover, given \(A,A'\in\mathcal{A}\), we made sure that only elements in \(M[G_{A_\cap A'}]\) have both an \(A\)-name and an \(A'\)-name, and therefore this model equals the intersection \(M[G_A]\cap M[G_{A'}]\).
		
		Finally, let us consider some \(B\subseteq I\) with \(B\notin\mathcal{A}\).
		By assumption there is a finite \(B'\subseteq B\) with \(B'\in\mathcal{B}\).
		Consider the reals \(r_i\) for \(i\in B'\).
		According to our construction we can decode bits of \(z\) by looking at elements of the intersection \(\bigcap_{i\in B'}r_i\) and whether \(r_{i_{B'}}\) or one of the other reals has the smaller next element.
		However, we should be careful that the coding points in \(\bigcap_{i\in B'} r_i\) were only added when intended.
		Inspecting the construction, the only time a point is added to \(\bigcap_{i\in J}r_i\) is if \(J\) is covered by an element of \(\mathcal{A}\) (and we met a dense set), covered by an element of \(\mathcal{B}\) (and we added an element to \(\bigcap_{i\in J}r_i\) on purpose), or if in the current step of the construction some of the coordinates in \(J\) occur for the first time (after a coordinate appears, the uniformity requirement guarantees that no new elements of the intersection are added accidentally).
		In the case of \(J=B'\in\mathcal{B}\), it clearly cannot be covered by any element of \(\mathcal{A}\) and the only element of \(\mathcal{B}\) that covers it is itself.
		Therefore the only time an element of \(\bigcap_{i\in B'}r_i\) can be added accidentally is if some coordinate in \(B'\) was mentioned for the first time in the given step.
		Since \(B'\) is finite there will only be finitely many such extraneous coding points and we will still be able to recover (a tail of) the real \(z\), ensuring nonamalgamability.
	\end{proof}
	
	An argument combining ideas from the above proof with those from the proof of \zcref{stat:CohenProjectionsNonamalg} shows that the same result can be obtained not just for filter-based Mathias forcing but also for families of posets that merely project to the same filter-based Mathias forcing.
	
	A careful reading of the proof of \zcref{stat:MathiasNonamalg} shows that the argument works in greater generality.
	Namely, we can allow the filter \(\F\) to vary between coordinates, as long as we keep this family of filters compatible.
	
	\begin{theorem}
		Suppose that a family \(\mathcal{A}\) is defined in \(M\) by a set of finite obstacles on a set \(I\).
		Let \(\set{\F_i}{i\in I}\in M\) be a family of filters on \(\omega\) extending the cofinite filter and linearly ordered under inclusion.%
		\footnote{Equivalently, the \(\F_i\) all extend to the same nonprincipal ultrafilter.}
		Then there are generic filters \(G_i\subseteq \M_{\F_i}\) over \(M\), for \(i\in I\), with the following properties:
		\begin{enumerate}
			\item If \(A\in\mathcal{A}\) then \(G_A\) is generic for \(\M_{\F_A}\) over \(M\).
			\item If \(B\in M\), \(B\subseteq I\), and \(B\notin\mathcal{A}\) then the family \(\set{M[G_i]}{i\in B}\) does not amalgamate in the generic multiverse.
			\item If \(A,A'\in\mathcal{A}\) then \(M[G_A]\cap M[G_{A'}]=M[G_{A\cap A'}]\).
		\end{enumerate}
		All products above are taken with finite support.
	\end{theorem}
	
	\begin{proof}
		The proof of \zcref{stat:MathiasNonamalg} transfers almost verbatim.
		We can no longer deal with uniform conditions, but that is a minor inconvenience.
		The only care needed is to see that we can still add coding points and that these don't get added too often.
		
		Given a condition \(p_n\) and a finite obstacle \(B\), we added a coding point by extending the stems of \(p_n\rest B\) by a common point.
		This is still possible in our current situation since the finitely many filters \(\F_i\) for \(i\in B\) are linearly ordered, which implies that the upper parts \(X_n(i)\) for \(i\in B\) have an infinite intersection (lying in the largest filter among the \(\F_i\)).
		Therefore we can not only find a point in the intersection (and additional ones to code the next bit of \(z\)), but we can find these coding points above \(\max\bigcup_{i\in B}s_n(i)\).
		Doing the coding this high will avoid adding undesired coding points with respect to other obstacles \(B'\).
	\end{proof}
	
	As an application we can construct nonamalgamable arrangements mixing both Cohen and filter-based Mathias forcing.
	We mention the simplest version of the argument, constructing a nonamalgamable pair of one Cohen real and one filter-based Mathias real.
	
	\begin{corollary}
		There are generic filters \(G\subseteq \Add(\omega,1)\) and \(H\subseteq \M_\F\) such that \(M[G]\) and \(M[H]\) do not amalgamate in the generic multiverse.
	\end{corollary}
	
	The key is that Cohen forcing is forcing equivalent to a specific filter-based Mathias forcing, the one where the filter is simply the cofinite filter on \(\omega\).
%
%
%
	
	\bibliographystyle{amsalpha}
	\bibliography{bibbase}

\newcommand{\etalchar}[1]{$^{#1}$}
\providecommand{\bysame}{\leavevmode\hbox to3em{\hrulefill}\thinspace}
\providecommand{\MR}{\relax\ifhmode\unskip\space\fi MR }
\providecommand{\MRhref}[2]{%
  \href{http://www.ams.org/mathscinet-getitem?mr=#1}{#2}
}
\providecommand{\href}[2]{#2}
\begin{thebibliography}{HHK{\etalchar{+}}19}

\bibitem[Abr10]{Abraham:Handbook}
Uri Abraham, \emph{Proper forcing}, Handbook of set theory. {V}ols. 1, 2, 3,
  Springer, Dordrecht, 2010, pp.~333--394. \MR{2768684}

\bibitem[BJ95]{BartoszynskiJudah:book}
Tomek Bartoszy\'nski and Haim Judah, \emph{Set theory}, A K Peters, Ltd.,
  Wellesley, MA, 1995, On the structure of the real line. \MR{1350295}

\bibitem[Cum10]{Cummings:Handbook}
James Cummings, \emph{Iterated forcing and elementary embeddings}, Handbook of
  set theory. {V}ols. 1, 2, 3, Springer, Dordrecht, 2010, pp.~775--883.
  \MR{2768691}

\bibitem[FGK19]{FriedmanGitmanKanovei2019:ModelOfSecondOrderArithmetic}
Sy-David Friedman, Victoria Gitman, and Vladimir Kanovei, \emph{A model of
  second-order arithmetic satisfying {AC} but not {DC}}, J. Math. Log.
  \textbf{19} (2019), no.~1, 1850013, 39. \MR{3960895}

\bibitem[FH09]{FuchsHamkins2009:DegreesOfRigiditySuslinTrees}
Gunter Fuchs and Joel~David Hamkins, \emph{Degrees of rigidity for {S}ouslin
  trees}, The Journal of Symbolic Logic \textbf{74} (2009), no.~2, 423--454.
  \MR{2518565}

\bibitem[Ham16]{Hamkins2016:UpwardClosureAndAmalgamationInGenericMultiverse}
Joel~David Hamkins, \emph{Upward closure and amalgamation in the generic
  multiverse of a countable model of set theory}, RIMS {Ky\^oky\^uroku} (2016),
  17--31.

\bibitem[HHK{\etalchar{+}}19]{HHKVW2019:SetTheoreticBlockchains}
Miha~E. Habi\v{c}, Joel~David Hamkins, Lukas~Daniel Klausner, Jonathan Verner,
  and Kameryn~J. Williams, \emph{Set-theoretic blockchains}, Archive for
  Mathematical Logic \textbf{58} (2019), no.~7-8, 965--997. \MR{4003645}

\bibitem[Jec03]{Jech2002:SetTheory}
Thomas Jech, \emph{Set theory}, millennium ed., Springer Monographs in
  Mathematics, Springer-Verlag, Berlin, 2003, The third millennium edition,
  revised and expanded. \MR{1940513}

\bibitem[JS91]{JudahShelah1991:ForcingMinimalDegree}
Haim Judah and Saharon Shelah, \emph{Forcing minimal degree of
  constructibility}, The Journal of Symbolic Logic \textbf{56} (1991), no.~3,
  769--782. \MR{1129141}

\bibitem[Mos76]{Mostowski1976:ModelsOfBGAxioms}
Andrzej Mostowski, \emph{A remark on models of the {G}\"{o}del-{B}ernays axioms
  for set theory}, Sets and classes (on the work by {P}aul {B}ernays) (Gert~H.\
  M\"{u}ller, ed.), Stud. Logic Found. Math., vol. Vol. 84, North-Holland,
  Amsterdam-New York-Oxford, 1976, pp.~325--340. \MR{446976}

\bibitem[Tod88]{Todorcevic1988:OscillationsOfReals}
Stevo Todor\v{c}evi\'c, \emph{Oscillations of real numbers}, Logic colloquium
  '86 ({H}ull, 1986), Stud. Logic Found. Math., vol. 124, North-Holland,
  Amsterdam, 1988, pp.~325--331. \MR{922115}

\end{thebibliography}
\end{document}